\documentclass[11pt,a4paper]{article}
\usepackage{ifthen,latexsym,amssymb,amsmath,bbm,fixmath,graphics}
\usepackage[nobysame,initials]{amsrefs}
\usepackage[shortlabels]{enumitem} 
\usepackage{url}
\usepackage{hyperref}
\usepackage{tikz}

\newcommand{\I}[1]{{\mathbbm #1}}

\renewcommand{\mid}{:}

\renewcommand{\ge}{\geqslant}
\renewcommand{\le}{\leqslant}

\newif\ifnotesw\noteswtrue

\newcommand{\hide}[1]{}

\newcommand{\beq}[1]{\begin{equation}\label{#1}}
\newcommand{\eeq}{\end{equation}}

\newtheorem{theorem}{Theorem}

\newcommand{\bpf}[1][Proof.]{\smallskip\noindent{\it #1.} }
\newcommand{\qed}{\nolinebreak\mbox{\hspace{5 true pt}%
  \rule[-0.85 true pt]{3.9 true pt}{8.1 true pt}}}
\newcommand{\epf}{\qed \medskip}

\newcommand{\ce}[4]{(#1,#2,#3,#4)}

\newcommand{\drawthreecube}[1]{
  \begin{tikzpicture}[scale=1.5]
    \foreach \a in {0,1} {
      \foreach \b in {0,1} {
        \foreach \c in {0,1} {
            \pgfmathsetmacro{\x}{\a + 0.5*\c}
            \pgfmathsetmacro{\y}{\b + 0.5*\c}
            \pgfmathtruncatemacro{\vnum}{4*\a + 2*\b + \c}
            \coordinate (V\vnum) at (\x,\y);
        }
      }
    }
    \foreach \a in {0,1} {
      \foreach \b in {0,1} {
        \foreach \c in {0,1} {
            \pgfmathtruncatemacro{\vnum}{4*\a + 2*\b + \c}
            \ifnum \a=0
              \pgfmathtruncatemacro{\vnumA}{4*1 + 2*\b + \c}
              \draw (V\vnum) -- (V\vnumA);
            \fi
            \ifnum \b=0
              \pgfmathtruncatemacro{\vnumB}{4*\a + 2*1 + \c}
              \draw (V\vnum) -- (V\vnumB);
            \fi
            \ifnum \c=0
              \pgfmathtruncatemacro{\vnumC}{4*\a + 2*\b + 1}
              \draw (V\vnum) -- (V\vnumC);
            \fi
        }
      }
    }
    \foreach \a in {0,1} {
      \foreach \b in {0,1} {
        \foreach \c in {0,1} {
            \pgfmathtruncatemacro{\vnum}{4*\a + 2*\b + \c}
            \filldraw[fill=white, draw=black] (V\vnum) circle (2pt);
        }
      }
    }
    \foreach \v in {#1} {
      \fill (V\v) circle (2pt);
    }
  \end{tikzpicture}
}

\newcommand{\drawfourcube}[1]{
  \begin{tikzpicture}[scale=1.2]
    \foreach \a in {0,1} {
      \foreach \b in {0,1} {
        \foreach \c in {0,1} {
          \foreach \d in {0,1} {
            \pgfmathsetmacro{\x}{\a + 0.5*\c + 1.8*\d}
            \pgfmathsetmacro{\y}{\b + 0.5*\c + 0.25*\d}
            \pgfmathtruncatemacro{\vnum}{8*\a + 4*\b + 2*\c + \d}
            \coordinate (V\vnum) at (\x,\y);
          }
        }
      }
    }
    \foreach \a in {0,1} {
      \foreach \b in {0,1} {
        \foreach \c in {0,1} {
          \foreach \d in {0,1} {
            \pgfmathtruncatemacro{\vnum}{8*\a + 4*\b + 2*\c + \d}
            \ifnum \a=0
              \pgfmathtruncatemacro{\vnumA}{8*1 + 4*\b + 2*\c + \d}
              \draw (V\vnum) -- (V\vnumA);
            \fi
            \ifnum \b=0
              \pgfmathtruncatemacro{\vnumB}{8*\a + 4*1 + 2*\c + \d}
              \draw (V\vnum) -- (V\vnumB);
            \fi
            \ifnum \c=0
              \pgfmathtruncatemacro{\vnumC}{8*\a + 4*\b + 2*1 + \d}
              \draw (V\vnum) -- (V\vnumC);
            \fi
            \ifnum \d=0
              \pgfmathtruncatemacro{\vnumD}{8*\a + 4*\b + 2*\c + 1}
              \draw[style=lightgray] (V\vnum) -- (V\vnumD);
            \fi
          }
        }
      }
    }
    \foreach \a in {0,1} {
      \foreach \b in {0,1} {
        \foreach \c in {0,1} {
          \foreach \d in {0,1} {
            \pgfmathtruncatemacro{\vnum}{8*\a + 4*\b + 2*\c + \d}
            \filldraw[fill=white, draw=black] (V\vnum) circle (2pt);
          }
        }
      }
    }
    \foreach \v in {#1} {
      \fill (V\v) circle (2pt);
    }
  \end{tikzpicture}
}

\author{Levente Bodn\'ar\thanks{Email: bodnalev@gmail.com
}\ \mbox{ and }Oleg Pikhurko\thanks{Email: O.Pikhurko@warwick.ac.uk}\\
Mathematics Institute and DIMAP\\
University of Warwick\\
Coventry CV4 7AL, UK}

\title{Some exact values of the inducibility and statistics constants for hypercubes}

\begin{document}

\maketitle

\begin{abstract}
We consider the following two problems in the standard hypercube graph on~$\{0,1\}^n$. The \emph{inducibility problem} for a given configuration $H\subseteq \{0,1\}^d$ (resp.\ the \emph{statistics problem} for given integers $d$ and $s$)
asks for the maximum, over all $S\subseteq\{0,1\}^n$, of the number of $d$-dimensional subcubes that intersect $S$ in a copy isomorphic to~$H$ (resp.\  contain exactly $s$ points from~$S$). Using flag algebras, we determine the limit of the corresponding maximal density of subcubes as $n\to\infty$ for 5 configurations $H\subseteq\{0,1\}^3$ and for 3 pairs $(d,s)$, namely for $(3,2)$, $(4,2)$ and $(4,4)$. Interestingly, the lower bounds in the last three cases come from blowups of small Hamming codes.

Our results on the inducibility problem were also independently proved by Rahil Baber.
\end{abstract}


\section{Introduction}\label{se:intro}

The \emph{$n$-hypercube} (or \emph{$n$-cube} for short) $Q_n$ is the graph on $\{0,1\}^n$, the set of binary sequences of length $n$, where two sequences are considered adjacent if they differ in exactly one coordinate. Thus $Q_n$ has $2^n$ vertices and $n\cdot 2^{n-1}$ edges. 

While $Q_n$ provides a very natural and widely applicable framework for studying set systems on $\{1,\dots,n\}$,  binary $n$-bit codes and Boolean $n$-variable functions, its graph structure is also of great interest on its own. Indeed, various versions of extremal graph problems were studied for hypercubes (see, for example,~\cite[Section~1]{GoldwasserHansen24} for references).

One direction of this type is to study possible distributions of a configuration $S\subseteq \{0,1\}^n$ among $d$-subcubes of~$Q_n$. Formally, a \emph{$d$-subcube} is a subset $C\subseteq \{0,1\}^n$ such that the subgraph $Q_n[C]$ of $Q_n$ induced by $C$ is isomorphic to $Q_d$. As it is easy to see, there are exactly ${n\choose d}\, 2^{n-d}$
$d$-subcubes and each can be obtained by fixing some $n-d$ coordinates and varying the $d$ remaining coordinates (called \emph{flip bits}) arbitrarily. We say that a $d$-subcube $C\subseteq \{0,1\}^n$ gives a \emph{copy} of $H\subseteq \{0,1\}^d$ in $S\subseteq \{0,1\}^n$ if there is a graph isomorphism $f$ from  $Q_d$ to $Q_n[C]$ with $f(H)=S\cap C$. In other words, copies of $H$ in $S$ are given by graph embeddings of $Q_d$ into $Q_n$ that send all vertices in $H$ into $S$ and all vertices in $\{0,1\}^d\setminus H$ into $\{0,1\}^n\setminus S$. Let $\Lambda(H,S)$ denote the number of copies of $H$ in $S$
and let the \emph{density} of $H$ in $S$ be
 $$
  \lambda(H,S):=\frac{\Lambda(H,S)}{\binom nd\, 2^{n-d}},
  $$
  which is the fraction of $d$-subcubes that give a copy of $H$ in $S$. For example, 
if $H=\{(0),(1)\}$ (and thus $d=1$) then, $\Lambda(H,S)$ is the number of edges of $Q_n$ that are inside~$S$.
  
These two parameters are the natural analogues of the subgraph count and the subgraph density for hypercubes. Motivated by a long and active line of research on the graph inducibility problem (see e.g.\ \cite{PippengerGolumbic75,BollobasNaraTachibana86,BrownSidorenko94,EvenzoharLinial15,HatamiHirstNorin14jctb,Hirst14,BaloghHuLidickyPfender16,KralNorinVolec19,Yuster19} for a small sample of results),
Goldwasser and Hansen~\cite{GoldwasserHansen21,GoldwasserHansen24} initiated the systematic study, for given $H\subseteq \{0,1\}^d$, of the \emph{inducibility function}
\beq{eq:Lambda}
 \Lambda(H,n):=\max\{\Lambda(H,S)\mid S\subseteq \{0,1\}^n\},\quad \mbox{for integer $n\ge d$},
 \eeq
 and the \emph{inducibility constant}
  $$
  \lambda(H):=\lim_{n\to\infty} \frac{\Lambda(H,n)}{\binom nd\, 2^{n-d}},
  $$
  where the limit exists since the ratio is non-increasing. (Indeed, for $n\ge d+1$, the density of $H$ in $S\subseteq \{0,1\}^n$ is the average of its densities inside $(n-1)$-dimensional subcubes.)

\begin{figure}
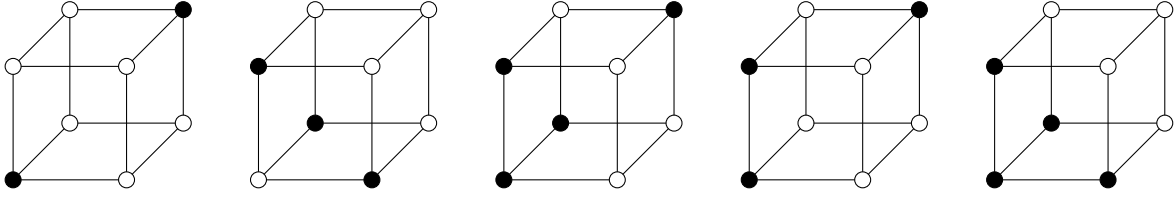

\begin{equation*}
    \drawthreecube{0,7} \qquad \drawthreecube{1,2,4} \qquad \drawthreecube{0,1,2,7} \qquad \drawthreecube{0,2,7} \qquad \drawthreecube{0,1,2,4}
\end{equation*}
    \caption{Configurations $W_7, W_8, W_9, W_{10}, W_{12}$ respectively.}\label{fi:1}
\end{figure}

Of course, $\Lambda(H,n)$ does not change if we replace $H\subseteq \{0,1\}^d$ by its image under an automorphism of $Q_d$ or by its complement $\{0,1\}^d\setminus H$. In this respect, there are $4$ non-equivalent configurations in $\{0,1\}^2$ and $14$ non-equivalent configurations in $\{0,1\}^3$. (We skip the trivial case $d=1$ when each inducibility constant is 1.) 
Goldwasser and Hansen~\cite{GoldwasserHansen21,GoldwasserHansen24} determined the inducibility constant (as well as the inducibility function in many cases) for 3 different configurations $H\subseteq \{0,1\}^2$ and 6 different configurations $H\subseteq \{0,1\}^3$. In the case $d=2$, the only open case is when $H=\{(0,0)\}$ consists of a single vertex of $Q_2$, with the best known lower bound being $\lambda(H)\ge 2/3$, see~\cite{GoldwasserHansen24}. As it was observed in~\cite{GoldwasserHansen24},  there are further configurations $H\subseteq \{0,1\}^3$ where the numerical upper bound on $\lambda(H)$ coming the flag algebra method
is close to a known lower bound. These are denoted by $W_7,W_8, W_9,W_{10},W_{12}$ in~\cite{GoldwasserHansen24} and are
\begin{eqnarray*}
 W_7&:=&\{(0,0,0),\,(1,1,1)\},\\ 
 W_8&:=&\{(0,0,1),\,(0,1,0),\,(1,0,0)\},\\
 W_9 &:=& \{(0,0,0),\, (0,0,1),\, (0,1,0),\, (1,1,1)\},\\
 W_{10} &:=& \{(0,0,0),\,(0,0,1),\, (1,1,1)\},\\
 W_{12} &:=& \{(0,0,0),\,(0,0,1),\,(0,1,0),\,(1,0,0)\},
 \end{eqnarray*}
 see also Figure~\ref{fi:1}.

Here, we show that the floating-point flag-algebra calculations can be converted into rational-valued matrices that prove the matching upper bound on the inducibility constant in each of these cases, in particular resolving Conjecture 9.1 (and Conjecture~5.3) from~\cite{GoldwasserHansen24}:

\begin{theorem}\label{th:W} We have the following  inducibility constants: $\lambda(W_7)=1/3$, $\lambda(W_8)=2/3$, $\lambda(W_9)=4/9$, $\lambda(W_{10})=5/12$ and $\lambda(W_{12})=1/2$.
\end{theorem}

Our proofs of the  upper bounds in the above theorem are computer-assisted. A modified version of the SageMath software \cite{sagemath}, developed by the first author, can be used by the reader for both generating and verifying flag algebra certificates, see Section~\ref{se:upper} for further details. As we learned after the appearance of our preprint, the upper bounds of Theorem~\ref{th:W} (as well as various upper bounds on the inducibility constant of different $H\subseteq\{0,1\}^4$) were proved by Rahil Baber back in 2013 (unpublished); his certificates and verifier code can be found in~\cite{Baber13data}.

Additionally, we were able to prove that $\lambda(W_{12})=1/2$ without using any computer calculations; this also follows from a more general result of John Goldwasser and Ryan Hansen (personal communication). In fact, we can determine the whole inducibility function for $W_{12}$:

\begin{theorem}\label{th:W12} We have $\Lambda(W_{12},3)=1$ and, for each $n\ge 4$, it holds that
	\beq{eq:W12}	
	\Lambda(W_{12},n)= \frac12\cdot {n\choose 3}\cdot 2^{n-3}.
	\eeq
\end{theorem}	

Thus, there are only three non-equivalent $H\subseteq \{0,1\}^3$ for which we do not know the inducibility constant:
 \begin{eqnarray*}
 W_3 &:=& \{(0,0,0)\},\\
 W_4 &:=& \{(0,0,0),\,(1,0,0)\}, \\ 
 W_5 &:=& \{(0,0,0),\,(0,0,1),\,(0,1,0)\}.
 \end{eqnarray*}
 The best known lower bounds are $\lambda(W_3)\ge 1/2$, $\lambda(W_4)\ge 4/9$ and $\lambda(W_5) \ge 1/3$ (see~\cite{GoldwasserHansen24}). However, the numerical bounds coming from flag algebra calculations do not match these. For $\lambda(W_5)$, the numerical upper bound is approximately $1/3 + 6.4 \times 10^{-5}$. While this is close to the lower bound of $1/3$, the difference does not appear to be just a floating-point rounding error. So closing this gap is likely to require new ideas.

\begin{figure}
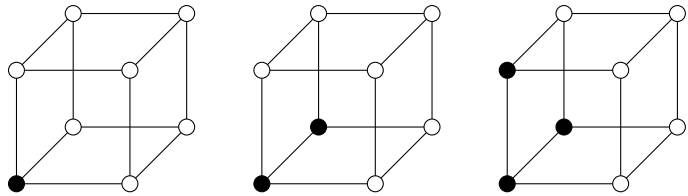

\begin{equation*}
    \drawthreecube{0} \qquad \drawthreecube{0,1} \qquad \drawthreecube{0,1,2}
\end{equation*}
    \caption{Unsolved configurations, $W_3, W_4$ and $W_5$.}\label{fi:2}
\end{figure} 
 
The related hypercube \emph{statistics problem} was recently introduced by Alon, Axenovich and Goldwasser~\cite{AlonAxenovichGoldwasser24}. It asks, for given integers $d$ and $s$ with $0\le s\le 2^d$, to determine the \emph{statistics function}
 \[
 \Lambda(n,d,s):=\max\{\Lambda(S,d,s)\mid S\subseteq \{0,1\}^n\},\quad \mbox{for integer $n\ge d$,}
 \]
 where $\Lambda(S,d,s)$ denotes the number of $d$-subcubes $C\subseteq \{0,1\}^n$ with $|S\cap C|=s$. In other words, we maximise, over all configurations $S\subseteq \{0,1\}^n$, the number of $d$-subcubes of $Q_n$ that contain exactly $s$ points from~$S$. Also, let the \emph{statistics constant}
 $$
 \lambda(d,s):=\lim_{n\to\infty} \frac{\Lambda(n,d,s)}{\binom nd\, 2^{n-d}}
 $$
 be the asymptotically maximal possible density of $d$-subcubes with exactly $s$ configuration points. Again, the limit exists as the ratio is non-increasing. 
 
It is trivial to see that $\lambda(d,0)=\lambda(d,2^{d-1})=1$ and, by passing to the complements, that $\lambda(d,s)=\lambda(d,2^d-s)$, so we can assume that $1\le s<2^{d-1}$. If $s=1$ then, up to automorphism, $Q_d$ contains only one possible single-point configuration and we get an instance of the inducibility problem (which, as we mentioned above, is open even for $d=2$). As observed in~\cite{Alon25talk}, further motivation for studying $\Lambda(n,d,1)$ comes from erasure codes: by a trivial reformulation of the definition of $\Lambda(n,d,1)$, we look for a code $S\subseteq\{0,1\}^n$ that maximises the number of pairs $(x,D)$, where $x$ is an element of $S$ and $D$ is a $d$-set of indices such that $x$ is uniquely determined even when we \emph{erase} each $x_i$ (i.e.,\ replace $x_i$ by a wildcard symbol `$*$') with $i\in D$.
 
While the initial paper~\cite{AlonAxenovichGoldwasser24}
mainly concentrated on various asymptotics as $d\to\infty$ (exhibiting rather different behaviour in many respects when compared to a related edge-statistic problem for graphs, see~\cite{KwanSudakovTran19,MMNT19,AlonHefetzKrivelevichTyomkyn20,FoxSauerman20}), computing the statistics function for small concrete $d,s$ is also of interest. In fact, Noga Alon~\cite{Alon25talk} posed the problem to determine $\lambda(d,s)$ for a single pair $(d,s)$ with $0<\lambda(d,s)<1$.  Here, we solve the following cases:

\begin{theorem}\label{th:statistics} It holds that $\lambda(3,2)=8/9$, $\lambda(4,2)=264/343$ and $\lambda(4,4)=26/27$.\end{theorem}

We use the same SageMath package to prove upper bounds (with the lower bounds coming from a special case of the construction appearing in the proof of~\cite[Theorem 2]{AlonAxenovichGoldwasser24}). Rahil Baber (personal communication) was able to re-prove Theorem~\ref{th:statistics} using his code.

\medskip\noindent\textbf{Organisation of the paper.} The rest of this paper is organised as follows.  In Section~\ref{se:lower}, we present the constructions from \cite{AlonAxenovichGoldwasser24,GoldwasserHansen24} that give the lower bounds in Theorems~\ref{th:W} and~\ref{th:statistics}. Section~\ref{se:upper} 
briefly discusses the flag algebra method for hypercubes and presents details on how to verify the stated upper bounds. Section~\ref{se:W12} contains the proof of Theorem~\ref{th:W12}. Finally, Section~\ref{se:conclude} contains some concluding remarks.

\section{Lower bounds}\label{se:lower}

First, we present a construction that is sufficiently general to cover all lower bounds in Theorems~\ref{th:W} and~\ref{th:statistics}.

For an integer $r\ge 1$, let $\I Z_{r}:=\{0,\dots,r-1\}$ denote the set of residues modulo $r$, while $[r]$ denotes the set $\{1,\dots,r\}$. Suppose that we have positive integers $r_1,\dots,r_t$ and a partition $[n]=A_1\cup \dots\cup A_t$. We define 
 \newcommand{\asum}[3]{\Sigma_{#1}^{
 #3}(#2)}
 \[
 \asum{A_1,\dots,A_t}{x_1,\dots,x_n}{r_1,\dots,r_t}:=\left(\sum_{j\in A_1} x_j\bmod{r_1},\dots,\sum_{j\in A_t} x_j \bmod{r_t}\right),\quad\mbox{for $(x_1,\dots,x_n)\in\{0,1\}^n$.}
 \]
 Thus, $\asum{A_1,\dots,A_t}{x_1,\dots,x_n}{r_1,\dots,r_t}$ is the element of $\I Z_{r_1}\times\cdots \times \I Z_{r_t}$ whose $i$-th component for $i\in [t]$ is  the sum modulo $r_i$ of those entries $x_j$ whose coordinate $j$ is in $A_i$. For a $d$-subcube $C\subseteq \{0,1\}^n$ with the $d$-set $D\subseteq [n]$ of flip coordinates, define
\hide{ $$
  \asum{A_1,\dots,A_t}{C}{r_1,\dots,r_t}:=\left(\sum_{i\in A_1\setminus D} x_i\bmod{r_1},\dots,\sum_{i\in A_t\setminus D} x_i \bmod{r_t}\right),
 $$
 for any $(x_1,\dots,x_n)\in C$.
} 
 $$ 
  \asum{A_1,\dots,A_t}{C}{r_1,\dots,r_t}:=\asum{A_1,\dots,A_t}{x_1,\dots,x_n}{r_1,\dots,r_t},
  $$
  where $(x_1,\dots,x_n)$ is the unique element of $C$ with $x_i=0$ for all $i\in D$. In other words, we just take the residues of the sums of the fixed coordinates of $C$ inside each part~$A_i$.
  
 The \emph{$(r_1,\dots,r_t)$-blowup} of a subset $T\subseteq \I Z_{r_1}\times\dots \times \I Z_{r_t}$ with respect to $(A_1,\dots,A_t)$ is the set 
 $$
 S:=\left\{x\in\{0,1\}^n\mid \asum{A_1,\dots,A_t}{x}{r_1,\dots,r_t}\in T\right\}.
 $$
 Informally speaking, we decide on the inclusion of a binary sequence into $S$ depending solely on the number of 1's  modulo~$r_i$ inside each group $A_i$, with $T$ denoting the set of those residue sequences that pass the test.

Note that blowups consist by definition of binary sequences only, even if some $r_i$'s are larger than~2.
For example, if $t=1$ (and thus $A_1=[n]$) then we include exactly those sequences whose \emph{weight} (the total number of 1's) taken modulo $r_1$ belongs to the given set $\{m\in \I Z_{r_1}\mid (m)\in T\}$ of residues.

Let us describe the constructions from~\cite{GoldwasserHansen24} that give the lower bounds in Theorem~\ref{th:W}. For brevity, we just give the values of $\underline{q} = (q_1, ..., q_t)$, $\underline{r} = (r_1, ..., r_t)$ and $T$, to represent the construction described above with part sizes $|A_i| = (q_i+o(1)) n$ for $i\in [t]$ as $n\to\infty$. 
 Since we prove a lower bound on $\lambda(W_i)$, we list only those $3$-subcubes that give a copy of $W_i$, without formally verifying that no other copies of $W_i$ are present. We start with the constructions that are easier to describe. 

\begin{itemize}
    \item For $W_7$ and $W_8$, take $\underline{q}:=(1), \ \underline{r}:=(3), \ T := \{(0)\}$ that is, we let  $S\subseteq \{0,1\}^n$ consist of those binary $n$-sequences whose weight is divisible by~$3$. 
    Then every $3$-subcube $C$ gives a copy of either $W_7$ or $W_8$, with $W_7$ occurring exactly when $\asum{[n]}{C}{3}=0$ (that is, the minimum weight of an element of $C$ is divisible by 3). Thus, $\lambda(W_7)\ge 1/3$ and $\lambda(W_8)\ge 2/3$.
    
    \item For $W_{12}$, take $\underline{q}:=(1), \ \underline{r}:=(4), \ T := \{(0), (1)\}$ that is, let  $S\subseteq \{0,1\}^n$ consist of those binary $n$-sequences whose weight modulo 4 is 0 or 1. Note that every $3$-subcube $C\subseteq \{0,1\}^n$ with $\asum{[n]}{C}{4}\in \{0,2\}$ gives a copy of $W_{12}$ (as then the intersection of $C$ with $S$ consists of exactly either its first or its last two layers). Thus $\lambda(W_{12})\ge 1/2$.

    \item For $W_9$, take $\underline{q}:=\left(\frac23, \frac13 \right), \ \underline{r}:=(4, 2), \ T := \{(0, 0), (1, 0), (2, 1), (3, 1)\}$. Thus, a binary sequence $(x_1,\dots,x_n)$ is added to $S$ if $(\sum_{i\in A_1}x_i\bmod 4,\sum_{i\in A_2}x_i\bmod2)$ belongs to $T$. If a $3$-subcube $C$ has exactly 2 flip coordinates in $A_1$ 
    then $C$ gives a copy of $W_9$, regardless of the value of
    $\asum{A_1,A_2}{C}{4,2}$. Thus $\lambda(W_9)\ge 4/9$.

    \item For $W_{10}$, take $\underline{q}:=\left( \frac16, ..., \frac16 \right), \ \underline{r}:=(2, ..., 2)$ and $T$ consisting of certain 24 elements of $\I Z_2^6$. We refer the reader to~\cite[Proposition~7.4]{GoldwasserHansen24} for details.
\end{itemize}

Next, we present constructions (from~\cite[Proof of Theorem~2]{AlonAxenovichGoldwasser24}) that give the lower bounds in our Theorem~\ref{th:statistics}, with the derivation of the lower bounds that they give following shortly.

\begin{itemize}
    \item For $\lambda(3, 2)$ and $\lambda(4, 4)$ take $\underline{q}:=\left( \frac13, \frac13, \frac13 \right), \ \underline{r}:=(2, 2, 2), \ T := \{(0, 0, 0), (1, 1, 1)\}$. In other words, a binary $n$-sequence is included into $S$ if the three parities of its sums inside $A_1$, $A_2$ and $A_3$ are all the same.
    \item For $\lambda(4, 2)$ take $\underline{q}:=\left( \frac17, ..., \frac17 \right), \ \underline{r}:=(2, ..., 2)$ and let $T\subseteq \I Z_2^7$ be the $(7, 4)$-Hamming code.
\end{itemize}

Since the construction $\underline{q}:=\left( \frac13, \frac13, \frac13 \right), \ \underline{r}:=(2, 2, 2), \ T := \{(0, 0, 0), (1, 1, 1)\}$ seems to give a very good lower bound on $\lambda(d,2^{d-2})$, let us calculate this more general density in the corresponding blowups $S\subseteq\{0,1\}^n$.  As before, we list only those subcubes that span a desired configuration.
Note that every $d$-cube $C\subseteq\{0,1\}^n$ such that its set $D$ of flip coordinates intersects at least two parts $A_i$ gives a copy of a $(d,2^{d-2})$-configuration in $S$. Indeed, take any two flip indices $i,j$ in two different parts, say $i\in A_1$ and $j\in A_2$, and note that, for any assignment of values to the remaining $d-2$ coordinates, there is the unique value of $(x_i,x_j)$ that gives an element of $C\cap S$: the value of $x_i$ (resp.\ $x_j$) has to be the unique binary bit that makes the parity inside $A_1$ (resp.\ $A_2$) equal to the parity of $A_3$. Thus $C$ contains exactly $2^{d-2}$ points inside $S$, as desired.
Since the probability that a random $d$-cube in $\{0,1\}^n$ has all flip coordinates inside one part is $3\cdot 3^{-d}$, we conclude that, for every integer $d\ge 2$, it holds that
 \beq{eq:d-2}
 \lambda(d,2^{d-2})\ge 1-\frac{1}{3^{d-1}}.
 \eeq
 
For $d=3$ and $4$, we obtain the lower bounds stated in Theorem~\ref{th:statistics}. For $d=2$, we get the best known lower bound $\lambda(2,1)\ge 2/3$. 
A different construction attaining it with $\underline{q}:=( 1 ), \ \underline{r}:=(3), \ T := \{(0)\}$, was presented in~\cite{GoldwasserHansen24}. 
\hide{On the other hand, Baber~\cite{Baber13data} obtained a computer-generated proof that $\lambda(2,1)\le 24/35$ (and this result can be reproduced using our package). It is an intriguing open question to close this gap.}
  

Since the construction for $\lambda(4, 2)$ seems to give rather good bounds on $\lambda(d,2^{d-3})$, we calculate this more general density. First, we need to recall the definition of the $(7,4)$-Hamming code $T$.  Let $N:=\I Z_2^3\setminus\{(0,\dots,0)\}$ be the set of non-zero vectors in the 3-dimensional vector space $\I Z_2^3$ over the field~$\I Z_2$. Clearly, $|N|=7$. We view $T$ as a subset of $\I Z_2^N$ and define it as the set of those functions $f:N\to\I Z_2$ such that the linear combination $\sum_{v\in N} f(v)\,v$ is the zero vector in~$\I Z_2^3$. 
It will be convenient to label the 7 sets partitioning $[n]$ by $N$; thus $[n]=\cup_{v\in N} A_v$. For each $i\in [n]$, let $v_i$ be the unique index $v\in N$ such that $i\in A_v$. 
Thus, by the definition of blowup, $(x_1,\dots,x_n)\in \{0,1\}^n$ is included into $S$ if and only $\sum_{i=1}^n x_i v_i$ is the zero vector in $\I Z_2^3$, where we view each $x_i\in\{0,1\}$ as an element of the field $\I Z_2$ (a scalar equal to 0 or 1). 

Consider a $d$-subcube $C$ with the set $D$ of active coordinates. Observe that if the vectors $v_i$ for $i\in D$ span the whole space $\I Z_2^3$ then $C$ gives a $(d,2^{d-3})$-configuration. Indeed, fix any 3 elements $i,j,h\in D$ with $v_i,v_j,v_h$ forming a basis of $\I Z_2^3$. As before, for any values of the remaining $d-3$ flip bits, there is the unique choice of scalars $(x_i,x_j,x_h)\in \{0,1\}^3$ that gives an element of $S$; thus $C$ contains exactly $2^{d-3}$ elements from $S$. 

If $n\to\infty$ and each part $A_i\subseteq[n]$ has $(\frac17+o(1))n$ elements (and $d$ is fixed), then the distribution of $(v_{i_1},\dots,v_{i_d})\in N^d$ for a random permutation of a random $d$-subset $\{v_{i_1},\dots,v_{i_d}\}$ of $[n]$ approaches the uniform distribution, where we choose each coordinate uniformly in $N$ independently of the others. In the limit, the probability of $d$ uniform elements of $N$ not containing a basis of $\I Z_2^3$ can be calculated as follows. The vector space $\I Z_2^3$ has $7$ planes (2-dimensional subspaces) and the probability that all $d$ vectors from $N$ land inside a given plane $P$ is $3^d/7^d$. We overcount each outcome when all $d$ vectors end on a line (in this case actually meaning that all $d$ vectors are the same) exactly 3 times (the number of different planes containing the line). Since there are $7$ lines, we have to subtract $7\cdot (3-1)\cdot 7^{-d}$ to correct this overcount.

We conclude that, for every $d\ge 3$, it holds that
 \beq{eq:d-3}
 \lambda(d,2^{d-3})\ge 1-\frac{3^d-2}{7^{d-1}}.
 \eeq
 For $d=4$, this gives the lower bound $\lambda(4,2)\ge 264/343$ stated in Theorem~\ref{th:statistics}.

\section{Upper bounds}\label{se:upper}




We use the flag algebra method of Razborov~\cite{Razborov07}, which is also described in e.g.~\cite{Razborov10,BaberTalbot11,SFS16,GilboaGlebovHefetzLinialMorgenstein22}. In the hypercube setting, we organise calculations in the same way as is done by
Baber~\cite[Section 2.1]{Baber12} (see also~\cite{Baber11thesis}) and Balogh, Hu, Lidicky and Liu~\cite{BaloghHuLidickyLiu14}; we refer the reader to~\cite{Baber12,BaloghHuLidickyLiu14} for the details. 

Very briefly, flag algebras can provide a proof system for working with limiting subcube densities in an unknown configuration $S\subseteq\{0,1\}^n$ as $n\to\infty$. In particular, an asymptotic inequality, such as an upper bound $\lambda(H)\le \lambda_0$, can be mathematically proved by representing the difference $\lambda_0-\lambda(H,S)$ as a sum of squares in the flag algebra framework (where this identity does not involve $n$ and has to hold asymptotically for every $S\subseteq \{0,1\}^n$ as $n\to\infty$). This amounts to solving a (typically rather large) semi-definite program, which can be numerically done on a computer. Then one needs to \emph{round} the found floating-point solution, that is, to find positive semi-definite matrices with rational coefficients that yield a flag algebra identity satisfying 
all required properties.%
\hide{Briefly, the computer generates all sum-of-squares relations up to a given size and finds a floating-point positive semi-definite matrix describing an usage of the sum-of-squares method to get the best bound. An optimal matrix is found using the CSDP  solver~\cite{Borchers01011999}, which then gets rounded to a rational positive semi-definite matrix, to provide an exact bound. 
}

This is fully automated in the modified version of SageMath (using the CSDP  solver~\cite{Borchers01011999}). While this package is still under development, a short guide on how to install it and its current functionality can be found in the GitHub repository \href{https://github.com/bodnalev/sage}{\url{https://github.com/bodnalev/sage}}. The calculations that were sufficient for the stated bounds used the `plain' flag algebra method (i.e.,\ without need for any extra assumptions coming from  some external considerations that reduce the set of possible candidates~$S$). 
The script used for this paper and the generated certificates can be found in the ancillary folder of the arXiv version of this paper or in a separate GitHub repository \href{https://github.com/bodnalev/supplementary_files}{\url{https://github.com/bodnalev/supplementary_files}}.

To get the correct upper bound for $W_7, W_9, W_{12}$, the multiplication tables (sum-of-squares relations) up to the 3-cube were enough, while the other results used the multiplication tables up to the 4-cube. The ancillary folder contains both binary and human-readable text certificates. The text files contain:
\begin{itemize}
    \item the final bound,
    \item the list of all base flags (that is, flags with the empty type),
    \item the list of all corresponding slack values,
    \item the list of typed flag vectors whose square provides the claimed bounds after averaging.
\end{itemize}

The provided script (\texttt{hypercube.ipynb}) can be used to generate all certificate files (which may be different from ours due to possibly different numeric outputs by SDP solvers). 
The script also contains a verifier, to check that the provided matrices with rational entries are indeed positive semi-definite and to calculate the exact bound they provide. Finally, the script provides functionality to interactively inspect the certificates, printing out a requested square used in the proof.

\section{Inducibility function for $W_{12}$}\label{se:W12}

Here we establish the stated lower bound on the inducibility function of $W_{12}$. Recall that this configuration consists of the elements of $\{0,1\}^3$ of weight at most 1.

\bpf[Proof of Theorem~\ref{th:W12}] As the case $n=3$ is trivial, assume that $n\ge 4$.

For the stated lower bound, we take the same configuration $S\subseteq \{0,1\}$ that we used for the inducibility constant: the $(4)$-blowup of $T:=\{(0),(1)\}$ modulo~$4$. Recall that a $3$-subcube $C$ gives a copy of $W_{12}$ if its minimum-weight element $x$ has even parity. Every $x\in \{0,1\}^n$ having $2i$ entries equal to $1$ appears ${n-2i\choose 3}$ times this way. Thus
$$
 \Lambda(W_{12},n)\ge \sum_{i=0}^{\lfloor (n-3)/2\rfloor} {n\choose 2i} {n-2i\choose 3}
 =  {n\choose 3}  \sum_{i=0}^{\lfloor (n-3)/2\rfloor}  {n-3\choose 2i} =  {n\choose 3}\cdot \frac12\cdot 2^{n-3},
$$ 
as desired.

Let us turn to proving the stated upper bound. We use induction on $n\ge 4$ with the base case 
requiring some case analysis.

In the case $n=4$, we have to derive a contradiction from assuming that for some $S\subseteq \{0,1\}^4$ there are at least $5$ copies of $W_{12}$, say on the 3-subcubes $C_1,\dots,C_5$. For each $i\in [5]$, let the $i$-th \emph{root} $s_i$ be the unique vertex in $S\cap C_i$ such that all its 3 neighbours inside $C_i$ belong to~$S$. Likewise, the $i$-th \emph{anti-root} $a_i$  be the vertex in $C_i\setminus S$ such that all its 3 neighbours in $C_i$ are outside of~$S$. Of course, the set of roots as a subset of $S$ is disjoint from the set of anti-roots. 

Let us show that no root can be adjacent to an anti-root in~$Q_4$. Suppose on the contrary that this is false. By symmetry, assume that $s_1=\ce0000$ and $a_2=\ce0001$. Every neighbour of $s_1$ inside $C_1$ has to be in $S$, thus the subcube $C_1$ must be $\ce***0$ (where we denote the flip bits by stars). Likewise, $C_2=\ce***1$. Thus we know $S$ on $C_1\cup C_4=\{0,1\}^4$ completely. Up to symmetry, every other $3$-subcube is of the form $\ce0***$ or $\ce1***$, with these containing respectively $3$ and $5$ elements  from $S$. Thus  there are no other copies of $W_{12}$ in $S$, contradicting our assumption that $\Lambda(W_{12},S)\ge 5$.

Let us show that no vertex can be a root more than twice. Suppose that this is false, assuming by symmetry that  $s_1=s_2=s_3=\ce0000$, $C_1=\ce ***0$, $C_2=\ce**0*$ and $C_3=\ce*0**$. This determines which sequences belong to $S$, except for the two sequences in $\ce*111$. In particular, 
all weight-2 and weight-3 sequences except possibly $\ce0111$ are in $\{0,1\}^4\setminus  S$. At least one of the 5 roots, say $b_4$, has to be different from $\ce0000$ (since a vertex of $Q_4$ is in four 3-subcubes in total). Since $b_4\in S$ has at least 3 neighbours in $S$, it can neither have weight 1 nor be in $\ce*111$, a contradiction.

Next, let us show that, in fact, $s_i\not= s_j$ for any distinct indices $i,j\in [5]$. Suppose that this is false, assuming by symmetry that $s_1=s_2=\ce0000$, $C_1=\ce ***0$ and $C_2=\ce**0*$. This determines which points belong to $S$ except for the elements in $\ce**11$. Since no vertex can be root with multiplicity more than 2, there are at least two distinct roots, say $s_3$ and $s_4$, different from~$\ce0000$. Each of them belongs to $S$ and has at least 3 neighbours in $S$. This rules out all weight-1 and weight-3 sequences since each has three weight-2 neighbours (of which at most one  can be in $S$). For a similar reason, $\ce1111$ cannot be a root. Excluding sequences that we know to be in $\{0,1\}^4\setminus S$, we are left with the only possible sequence $\ce0011$ for $s_3$ and $s_4$, a contradiction.

Since $W_{12}$ is isomorphic to its complement $\{0,1\}^3\setminus W_{12}$, the above statement also applies to anti-roots, that is, no two anti-roots are the same.

By symmetry assume that $C_1=\ce***0$ and that $S\cap C_1$ consists exactly of sequences of weight at most 1. Then every weight-1 vertex $x$ in $C_1$ has at least two weight-2 neighbours in $C_1$ (which are outside of $S$). Thus $x$ cannot be a root and $C_1$ contains at most one root. Likewise $C_1$ contains at most one anti-root.  Thus the 3-subcube $\ce***1=\{0,1\}^3\setminus C_1$ has to contain at least 4 roots and at least 4 anti-roots, all distinct and thus exhausting the whole subcube. This forces an adjacency between these two groups, contradicting our earlier conclusion. This contradiction proves $\Lambda(W_{12},4)\le 4$, that is, the case $n=4$ of~\eqref{eq:W12}.

Suppose that we have already proved that~\eqref{eq:W12} holds for some $n\ge 4$. Take any configuration $S\subseteq \{0,1\}^{n+1}$. For each partition of $Q_{n+1}$ into two $n$-cubes, we have that at most $1/2$ of 3-cubes inside each part induce a copy of $W_{12}$. By averaging over all partitions, we conclude that at most $1/2$ of all 3-cubes inside $\{0,1\}^{n+1}$ induce a copy of $W_{12}$, as desired.\epf

\section{Concluding Remarks}\label{se:conclude}

It would be interesting to see if some of the above flag algebras certificates can be used to obtain not only the value of the inducibility or statistics constant but also the exact value of the corresponding function, at least for all sufficiently large~$n$. For a wide class of extremal graph problems, Pikhurko, Slia\v can and Tyros~\cite{PikhurkoSliacanTyros19} presented some sufficient conditions to proving such a result in an automated way and successfully applied it some instances of graph inducibility. The approach from~\cite{PikhurkoSliacanTyros19} crucially relies on an intermediate step of proving the \emph{stability property} which, roughly speaking, states that all almost extremal graphs have very similar structure. However, the considered hypercube problems, unless trivial, cannot satisfy the naive definition of stability: for example, we can cover almost all of $Q_n$ by a few vertex sets $V_i$ with every two at the graph distance more than $d$ from each other, and then replace an optimal configuration $S\subseteq \{0,1\}^n$ by letting the new configuration $S'$ on each $V_i$ be the intersection of $f_i(S)\cap V_i$ for some automorphism $f_i$ of~$Q_n$. Then $S'$ is almost optimal but, in general, it can be far away from any extremal configuration. A weaker form of stability, called the \emph{local stability} was proposed by Johnson and Talbot~\cite{JohnsonTalbot10}, who established it for a hypercube Tur\'an-type problem (\cite[Theorem 2]{JohnsonTalbot10}). While it is plausible that flag algebra proofs may be used to derive the local stability in some instances, it is unclear if this property can help with determining the extremal function exactly.

In some cases, the exact result can be directly derived from the asymptotic value by the blowup trick: for example, if the upper bound in~\eqref{eq:W12} fails for for some $n\ge 4$, that is, some $S\subseteq \{0,1\}^n$ has strictly more than $\frac12\cdot {n\choose 3}\cdot 2^{n-3}$ copies of $W_{12}$, then uniform blowups of $S$ would violate $\lambda(W_{12})\le 1/2$. However, this reduction, while giving some concrete bounds on the inducibility or statistics function, did not lead to the exact values in any of the other cases solved in the paper.

The configuration $T=\{(0,0,0),\,(1,1,1)\}$ that gives the lower bounds in~\eqref{eq:d-2} can be viewed as the $(3,1)$-Hamming code. Given the matching upper bounds of Theorem~\ref{th:statistics}, blowups of  the $(2^k-1,2^k-k-1)$-Hamming code $T$ seem good candidates to consider for the $(d,2^{d-k})$-statistics problem. 
However, the corresponding bounds are not always sharp. For example, the lower bound $\lambda(d,1)\ge \prod_{i=1}^{d-1}\frac{2^d-2^i}{2^d-1}$
coming from the uniform $(2,\dots,2)$-blowups of the $(2^d-1,2^d-1-d)$-Hamming code 
is not sharp for any $d\ge 2$:  if $d\le 4$ (resp.\ $d\ge 5$) then $(d+1)$-blowups of $T=\{(0)\}$ (resp.\ binomial $1/2^d$-random subsets of $\{0,1\}^n$) give the strictly better lower bound of $2/(d+1)$ (resp.\ $(1-1/2^d)^{2^d-1}$).

\begin{figure}
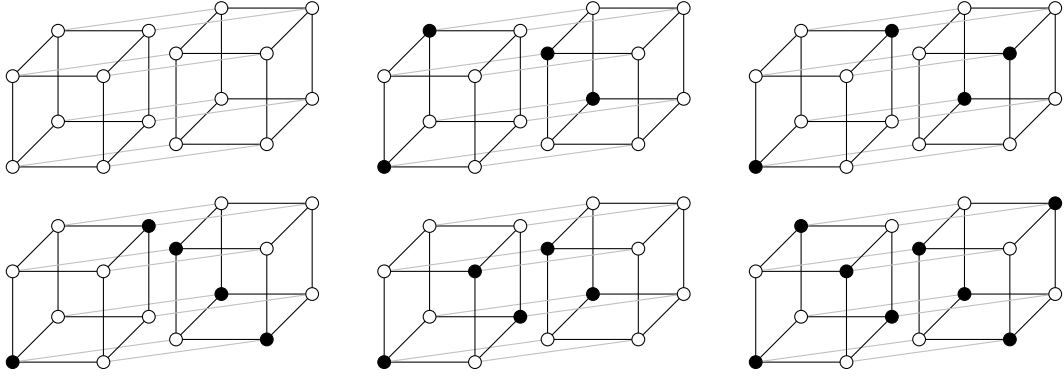

\begin{equation*}
\begin{split}
    \drawfourcube{} \qquad \drawfourcube{0,3,5,6} \qquad \drawfourcube{0,3,13,14} \\
    \drawfourcube{0,3,5,9,14} \qquad \drawfourcube{0,3,5,10,12} \qquad \drawfourcube{0,3,5,6,9,10,12,15}
\end{split}
\end{equation*}
    \caption{The configurations with $c_H=0$. The top row shows $H_1, H_2, H_3$, the bottom row shows $H_4, H_5, H_6$.}\label{fi:3}
\end{figure}

Determining the value of $\lambda(2,1)$ is an intriguing open problem. Rahil Baber~\cite{Baber13data} proved via flag algebras that $\lambda(2,1)\le 24/35$ (and this result can be reproduced using our package). While $24/35$ is the value of the flag algebra relaxation when using 4-cubes, we can show that it cannot be equal to $\lambda(2,1)$. In brief, the argument goes as follows. The proved upper bound of $24/35$ is derived from the following flag algebras identity:
 $$
 \frac{24}{35}-Z_2=\mathrm{SOS}+\sum_{H\subseteq \{0,1\}^4} c_H H,
 $$
 where $Z_2\subseteq \{0,1\}^2$ denotes the one-point configuration (whose density we bound from above), the  SOS (sum-of-squares) term is always non-negative, the last sum is over non-isomorphic $H$,
  and every coefficient $c_H$ is non-negative. It happens that there are exactly 6 non-isomorphic $H\subseteq \{0,1\}^4$ with $c_H=0$. See Figure \ref{fi:3} for these configurations. It follows that if a hypothetical configuration $S\subseteq\{0,1\}^n$ has $Z_2$-density $24/35+o(1)$ as $n\to\infty$ then $H_1,\dots,H_6$ are the only 4-dimensional configurations that can appear in $S$ with density not converging to 0. 
A computer can be used to test that all 5-dimensional extension of $H_5$ contains a 4-dimensional sub-configuration not from $H_1, \dots, H_6$. Thus $H_5$ must also have density $o(1)$ in~$S$. Now, if we run the flag algebra calculation of maximizing the density of $Z_2$ where we additionally require that the density of $H_5$ is 0 in the limit, then the returned numerical solution can be rounded to give an upper bound strictly less than $24/35$. The new bound applies to $S$, giving the desired contradiction. As we already mentioned, the best known lower bound is $\lambda(2,1)\ge 2/3$ (attained by two different constructions: $(3)$-blowup of, say  $T=\{(0)\}$, and the uniform $(2,2,2)$-blowup of the $(3,1)$-Hamming code $\{(0,0,0),(1,1,1)\}$).


\hide{
Also, let us observe that the inductive argument used in Theorem~\ref{th:W12} cannot be used to determine the statistics function $\Lambda(n,3,1)$
exactly, because it is strictly larger than $\frac12\cdot {n\choose 3}\cdot 2^{n-3}$ for every $n\ge 3$. Indeed, suppose that this fails for some $n$ and consider the $(4)$-blowups $S_i\subseteq \{0,1\}^n$ of $T=\{(i)\}$ for all four possible $i\in \I Z_4$. 
Then each $3$-subcube $C$ gives a copy of $(3,1)$ in exactly two configurations $S_i$: namely, if $i-\asum{[n]}{4}{C}$ is equal to 0 or 1 modulo 4. Thus the density of $(3,1)$-configurations in each $S_i$ is exactly $1/2$. 
It follows that each $S_i$ has exactly $2^{n-2}$ elements. A quick way to obtain a contradiction now is to observe that then the $(4)$-blowups of, say, $T=\{(0)\}$ would give that $\lambda(n,2^{n-2})=1$, contradicting~\cite[Theorem~1]{AlonAxenovichGoldwasser24}.
}

\newcommand{\DSRow}[3]{$#1$ &#2\\}
\begin{table}[t]
\begin{center}
    \begin{tabular}{ c |c}
         \DSRow{(d,s)}{Bound}{High precision floating point bound}
         \hline
         \DSRow{(3,1)}{0.6100426...}{0.61004260925419325294629225444}
         \DSRow{(3,3)}{0.6842325...} {0.68423255337451758617700846560}
         \DSRow{(4,1)}{0.6025329...}{0.60253292093316100763537868050}
         \DSRow{(4,3)}{0.6812635...}{0.68126354057902654997950548331}
         \DSRow{(4,5)}{0.7096421...}{0.70964218689054657579484995029}
         \DSRow{(4,6)}{0.8540599...}{0.85405998941874059355536504411}
         \DSRow{(4,7)}{0.7270395...}{0.72703955379766565090367509998}
    \end{tabular}
\end{center}
\caption{
Upper bounds on $\lambda(d,s)$ coming from flag algebras.}\label{ta:1}
\end{table}
Table~\ref{ta:1} lists the numerical upper bounds on $\lambda(d,s)$ suggested by flag algebra calculations for those pairs $(d,s)$ with $d\in \{3, 4\}$ and $s<2^{d-1}$ for which the exact value
is still unknown. 

As we already mentioned earlier, we made available not only the code that verifies our flag algebra certificates but also the code that we used to generate them. It should be fairly easy for the reader to modify the latter in order to do flag algebra calculations for some other extremal hypercube problems. We hope that the code will be useful for proving various further results.

\section*{Acknowledgements}

The authors thank Rahil Baber and John Goldwasser for useful comments, and  Rahil Baber for sharing his computer code and data.

Both authors were supported 
by ERC Advanced Grant 101020255.

\hide{
\section*{Various remaining results/observations (not for the paper)}

\hide{
This mathematica code indicates that modulo-4 constructions give a lower bound on $\Lambda(n,W_3)$ which is always greater than $\frac12 {n\choose d}2^{n-d}$ so we cannot resolve it as for $W_{12}$.
\begin{verbatim}
g[r_] := Binomial[n, d]*(Binomial[n - d, r] + Binomial[n - d, r - d])
d = 3; Table[{n, 
  Max[Table[
     Sum[cubesInBlowup[4 r + i], {r, 0, Floor[(n - i)/4]}], {i, 0, 
      3}]] - Binomial[n, d]*2^(n - d)/2}, {n, 4, 100}]
\end{verbatim}
}

From Noga Alon's talk (not in their paper)

\verb=https://player.vimeo.com/video/1055655586=

Open:

We know $\lambda(d,s)=1-\Theta(1/s)$ if $s$ is a power of $2$. What is the best possible constant in front of $1/s$? (The bounds in the paper give between 1/4 and 1 for $s=o(d)$ with $1/4$ decreasing when $o(d)<s\le d-o(d)$ so this question is a bit vague...)
}


\bibliography{bibexport}

\end{document}